\documentclass{article}

\usepackage{amsmath}
\usepackage{amsfonts}
\usepackage{amssymb}
\usepackage{color}  
\usepackage{amsthm}
\usepackage{graphicx}
\usepackage{caption}
\usepackage{comment}
\usepackage{subcaption}
\usepackage{algpseudocode}
\usepackage{algorithm}
\usepackage[utf8]{inputenc}
\usepackage{tabularx}
\usepackage{multirow}
\usepackage{hyperref}
\usepackage[english]{babel}
\usepackage{setspace}

\usepackage{natbib}
\bibpunct[, ]{(}{)}{,}{a}{}{,}%
\title{Graph Generation: A New Approach to Solving Expanded Linear Programming Relaxations}

\author{Julian Yarkony \textsuperscript{\rm 1, \rm 3}, Naveed Haghani\textsuperscript{ \rm2}, Amelia Regan\textsuperscript{ \rm 3}   \\[2ex] 

\textsuperscript{\rm 1} Laminaar Optimization Research Group, La Jolla, CA\\ 

\textsuperscript{\rm 2}University of Maryland, College Park, MD\\

\textsuperscript{\rm 3}University of California, Irvine, CA\\

}\date{October 2021}

\usepackage{natbib}
\usepackage{graphicx}

\begin{document}

\maketitle

\begin{abstract}
    In this article we introduce Graph Generation, an enhanced Column Generation (CG) algorithm for solving expanded linear programming relaxations of mixed integer linear programs. 
    To apply Graph Generation, we must be able to map any given column to a small directed acyclic graph for which any path from source to sink describes a feasible column. This structure is easily satisfied for vehicle routing and crew scheduling problems; and other such problems where pricing is a resource constrained shortest path problem. Such graphs are then added to the restricted master problem (RMP) when the corresponding column is generated during pricing. The use of Graph Generation does not weaken the linear programming relaxation being solved. At any given iteration of CG enhanced by Graph Generation; the technique permits the RMP to express a much wider set of columns than those generated during pricing, leading to faster convergence of CG. Graph Generation does not change the structure of the CG pricing problem. We show how the method can be applied in a general way, and then demonstrate the effectiveness of our approach on the classical Capacitated Vehicle Routing Problem.
    
\end{abstract} 
\section{Introduction}
Expanded linear programming (LP) relaxations of mixed integer linear programs provide much tighter relaxations than typical compact formulations for large classes of important problems in logistics  \citep{Desaulniers2005}, Chapter 1, and recently in computer vision/machine learning \citep{yarkony2020data,FlexDOIArticle}. Such formulations can often be described as follows: \par We are given a set of agents that must cover a set of tasks. Each agent completes a subset of the tasks called an assignment. The problem is to provide each agent with an assignment so as to ensure that each task is covered at least once and the total cost of the selected assignments is minimized. \par We should note that expanded LP relaxations efficiently permit additional constraints, which are often referred to as side constraints. The description of what constitutes a feasible assignment for a given agent or how to map that assignment to a cost is obscured in this formulation. However the revised simplex algorithm can still be used to solve this problem if for any given agent we can compute the lowest reduced cost assignment in a step called pricing. \par
In applications in logistics, the pricing operation is typically a resource constrained shortest path problem \citep{costa2019,Desrochers1992} or a knapsack problem \citep{diaz2002branch}. These are NP-hard \citep{karp_old} problems that can be solved efficiently in practice for problems of relevant scale \citep{barnprice}. In computer vision and machine learning these problems often correspond to solving a quadratic unconstrained binary optimization problem \citep{yarkony2020data}. It is frequently the case that these problems have special structures permitting fast solving in practice \citep{wang2018accelerating}.\par Expanded LP relaxations are solved using Column Generation (CG) \citep{barnprice}. CG solves expanded LP relaxations by mimicking the revised simplex. CG proceeds by solving the optimization over a subset of the assignments (primal variables) followed by generating assignments via pricing. The LP relaxation over this finite set is referred to as the restricted master problem (RMP).  Expanded LP relaxations are tighter than their corresponding compact formulation when the CG pricing problem has non-integer super-optimal solutions, that are feasible but for being non-integer \citep{desrosiers2005primer}.  

A primal variable generated in the course of CG is referred to as a column as it is associated with a new column in the constraint matrix. Each column is typically an assignment (a route for example). CG can suffer from slow convergence in problems where the number of non-zero entries in the constraint matrix commonly associated with a generated column exceeds a threshold number (for example, 8-12 stops in a route, as discussed in \citep{elhallaoui2005dynamic}). This difficulty can be circumvented with strong dual stabilization \citep{FlexDOIArticle,du1999stabilized}. \citep{haghani2020smooth} recently introduced Smooth dual optimal inequalities (S-DOI) for problems embedded in metric spaces such as vehicle routing and facility location problems. S-DOI permit low cost swap operations between nearby items in the primal. In the dual, S-DOI enforce that dual variables must change smoothly across space. This drastically decreases the size of the dual space that CG needs to search over, hence providing very large speedups.  S-DOI are expanded on in \citep{yarkony_Detour_DOI} to create Detour-dual optimal inequalities (DT-DOI). DT-DOI permit low cost swap operations between items in the route with items near to other items on that route. This formulation in essence maps each column to a cone of feasible columns (with increased costs beyond their true cost). This cone is then added to the RMP using a small set of additional primal variables/primal constraints; not by adding one primal variable for each column in the cone. The use of DT-DOI does not loosen the LP relaxation or alter the structure of the pricing problem. Thus instead of generating a single column, the pricing problem of CG generates a new set of variables/constraints for the RMP for which any solution over those variables describes a non-negative combination of columns. 

Detour-DOI are very powerful; but they are only able to be applied to a subset of problems in logistics associated with cost terms embedded on a metric space. Here we seek to generalize DT-DOI to add a cone of columns to the RMP in each iteration. This is a promising avenue for problems with varying structures. Such cones need to trade off the ability to express diverse sets of low cost columns not in the current RMP, with the need to not add large numbers of additional primal variables/constraints to the RMP. \par In this article we adapt this approach of adding cones of solutions to the RMP to arbitrary problems. Specifically we map each column produced during pricing to a small directed acyclic graph on which every path from source to sink corresponds to a feasible column with identical cost to the actual column. The larger this graph, the more columns can be expressed, but the more difficult computation becomes for the RMP. However in many applications, pricing rather than the RMP is the computational bottleneck \citep{desrosiers2005primer,haghani2020smooth,FlexDOIArticle} motivating the generation of better dual solutions to accelerate CG. We refer to our approach as \textbf {\large{Graph Generation}} (GG).  

We organize this document as follows. In Section \ref{Sec_litRev} we review the relevant literature on stabilizing CG. In Section \ref{Sec_basicCG} we provide formal description of standard CG. In Section \ref{Sec_GG} we introduce GG in an application agnostic manner. In Section \ref{Sec_CVRP_apply} we apply GG to the  classic Capacitated Vehicle Routing Problem (CVRP).  In Section \ref{Sec_exper} we provide experimental validation of the effectiveness of GG on CVRP. In Section \ref{Sec_conc} we conclude and discuss extensions. Before we begin we present the acronyms used most often in this paper in Table \ref{myTabNotatoin}. 

\begingroup
\setlength{\tabcolsep}{10pt} 
\renewcommand{\arraystretch}{1.5} 
\begin {table}[H]
\begin{center}
\begin{tabular}{|c|c|} 

 \hline
 CVRP & Capacitated Vehicle Routing Problem \\
 CG & Column Generation \\ 
 RMP & Restricted Master Problem \\
 MP & Master Problem \\
 GG & Graph Generation \\
 DOI & Dual Optimal Inequalities \\
 S-DOI & Smooth Dual Optimal Inequalities \\
 DT-DOI & Detour Dual Optimal Inequalities \\
 
 \hline
\end{tabular}
\end{center}
\caption {\bf{Acronyms}} \label{tab:title} 
\label{myTabNotatoin}
\end{table}
\endgroup
\section{Literature Review}
\label{Sec_litRev}

We now consider some approaches for dual stabilization of Column Generation (CG) related to our Graph Generation (GG) approach.  

\subsection{General Stabilization Methods}
  
Due to inherent instability of the dual variables in CG, many methods of stabilization have been proposed over the last two decades. Du Merle et al formalized the idea of stabilized CG in their 1999 paper of that name \citep{du1999stabilized}. That paper proposed a 3-piecewise linear penalty function to stabilize CG. Ben Amor and Desrosiers later proposed a 5-piecewise linear penalty function for improved stabilization \citep{amor2006proximal}. Shortly after, Oukil et al use the same framework to attack highly degenerate instances of multiple-depot vehicle scheduling problems \citep{oukil2007stabilized}. Ben Amor et al later proposed a general framework for stabilized CG algorithms in which a stability center is chosen as an estimate of optimal solution of the dual formulation \citep{amor2009choice}. Gonzio et al proposed a primal-dual CG method in which the sub-optimal solutions of the restricted master problem (RMP) are obtained using an interior point solver that is proposed in an earlier paper by the first author \citep{gondzio1995hopdm}. They examine their solution method relative to standard CG and analytic center cutting plane method proposed by Babonneau et al. \citep{babonneau2006solving, babonneau2007proximal}. They found that while standard CG is efficient for small problem instances, that the primal-dual CG method achieved the best solutions on larger problems \citep{Gondzio2013New}.

\subsection{Trust Region Based Stabilization}
Trust region based methods exploit the understanding of CG operating a search algorithm over the dual space \citep{marsten1975boxstep}. Such methods seek to maximize the Lagrangian Relaxation of the master problem.  Since columns generated at a given dual solution do not provide good information regarding the Lagrangian relaxation at distant points in the dual space CG can have a problem of bouncing between distant points in the dual space. This is circumvented by establishing a trust region computed at the point corresponding to the greatest Lagrangian relaxation thus far identified, and limiting search around it \citep{marsten1975boxstep}. This approach is extend in \citep{du1999stabilized} in which dual variables are penalized for leaving the trust region. 
%
\subsection{Dual Optimal Inequalities}

Dual optimal inequalities \citep{ben2006dual} provide provable bounds on the space where the  optimal dual solution lies. In this manner it reduces the size of the dual space that CG must search over. In the primal form these correspond to slack variables that may provide for swap operations between items \citep{haghani2020smooth} or provide rewards for over-covering items \citep{FlexDOIArticle}.  The corresponding primal variables are provably inactive in an optimal solution to the RMP.  

\subsection{Jigsaw Pricing}
Often it is desired to generate a set of columns that ``fit well together" in terms of producing integer solutions. This can also accelerate the convergence of CG. For a given dual solution generated by the RMP this is achieved  as follows. \\ We alternate between \textbf{(1)} solving pricing given the current dual variables; \textbf{(2)} setting to zero the dual value reward terms corresponding to covering items that have been covered by that column. We terminate this round of pricing when no column has negative reduced cost. Many variants of this have been produced, notably by exploiting the dynamic programming based solution approaches to pricing to generated multiple distinct columns. Our understanding is that this method is a heuristic employed by many practitioners of CG, though to date it has not been given a formal name in the literature. This concept was introduced to the first author via personal communication with Professor Jacques Desrosiers of the University of Montreal in early 2021. 

\subsection{Family Column Generation}

Family Column Generation \citep{haghani2021family}  builds off trust region based approaches. In Family Column Generation each column in the RMP is associated with a family of columns (often of exponential size in the number of master problem constraints) over which it is easy to price over.  Given the trust region each family corresponding to a column generated thus far is mapped to a column.  This column has the lowest possible reduced cost evaluated at an extreme point in the trust region (of columns inside that family) and has no greater reduced cost than the original column (column corresponding to the family) over all points in the trust region. The RMP is then solved over this trust region. This accelerates CG by ensuring that columns generated previously can play a role in stabilizing CG even if the trust region would normally ensure that their corresponding dual constraint is satisfied.  

\subsection{Detour-Dual Optimal Inequalities}
Detour-dual optimal inequalities (DT-DOI) \citep{yarkony_Detour_DOI} are applied to accelerate CG in problems where cost terms are embedded on a metric space such as the Capacitated Vehicle Routing Problem.  In the primal they permit swap operations between customers in the column and customers nearby any customer in that column. Such operations are expressed by making a detour. 
%
%
DT-DOI create a set of unique variables and constraints associated with that specific column upon generation.  These can be understood as providing specifications for a cone of columns with increased cost from their actual cost.  This increase is a function of the swap operations used to produce the column.  

\section{Review of Standard Column Generation}
\label{Sec_basicCG}
In this section we provide a formal mathematical review of the standard Column Generation (CG) approach for solving expanded linear programming (LP) relaxations. We describe expanded LP relaxations using the following notation. 

We use $\Omega$ to denote the set of primal variables, which we index by $l$. We use $A_{:l}$ to denote the column in the constraint matrix $A$ corresponding to primal variable $l$.  We use $\vec{b}$ to denote the vector of right hand side constants on the constraints. We express optimization using inequality constraints only since equality constraints can be rewritten as two inequality constraints. We use $\pi$ to denote the dual variable vector associated with the inequality constraints. We now frame the standard LP formulation with dual variables $\pi$ written next to the associated constraint in $[]$.  
\begin{subequations}
\label{orig_LP}
\begin{align}
    \min_{\theta \geq 0} \sum_{l \in \Omega }c_l\theta_l\\
    \sum_{l \in \Omega}A_{:l}\theta_l \geq \vec{b} \quad [\pi]
\end{align}
\end{subequations}
We refer to \eqref{orig_LP} as the master problem (MP).  Solving \eqref{orig_LP} is done using CG when $\Omega$ can not be explicitly enumerated (as is the case in many applications in operations research \citep{barnprice,lubbecke2005selected}); and the lowest reduced cost primal variable can be computed given any dual solution $\pi$. The computation of the lowest reduced cost column $l\in \Omega$ is written as follows using $\bar{c}_l$ to denote the reduced cost of column $l$.
\begin{subequations}
\label{pricing}
\begin{align}
    \min_{l \in \Omega}\bar{c_l}\\
    \bar{c}_l=c_l-\pi^{\top}A_{:l} \quad \forall l \in \Omega
\end{align}
\end{subequations}
%
Since $\Omega$ can not be explicitly enumerated, CG relies on solving \eqref{orig_LP} over a nascent subset of the variables $\Omega_R$. The aim is to construct $\Omega_R$ in a manner s.t. solving \eqref{orig_LP} over $\Omega_R$ provides the optimal solution to optimization as if we had considered all of $\Omega$.  CG solves \eqref{orig_LP} by alternating between \textbf{(a)} solving \eqref{orig_LP} over $\Omega_R$, which is referred to as the restricted master problem (RMP) and \textbf{(b)} computing the lowest reduced cost column using \eqref{pricing} and adding it to $\Omega_R$.  Termination of CG occurs when no negative reduced cost columns exist. When no negative reduced cost columns exist then the optimal solution provided by the RMP is optimal for the MP. CG is initialized with the columns composing any feasible solution to \eqref{orig_LP}. It can alternatively be initialized by artificial variables  to ensure a feasible solution exists in each iteration of CG.  These artificial variables are associated with prohibitively high cost to use so that they are not used in an optimal solution to \eqref{orig_LP}.  We now write the RMP over $\Omega_R$, which we define the solution to as $\Psi(\Omega_R)$. 
\begin{subequations}
\label{orig_RMP}
\begin{align}
    \Psi(\Omega_R)=\min_{\theta \geq 0} \sum_{l \in \Omega_R }c_l\theta_l\\
    \sum_{l \in \Omega_R}A_{:l}\theta_l \geq \vec{b} \quad [\pi]
\end{align}
\end{subequations}
In Alg \ref{basicCG} we provide the CG solution to \eqref{orig_LP} with annotation below.    
\begin{algorithm}[!b]
 \caption{Basic Column Generation}
\begin{algorithmic}[1] 
\State $\Omega_R\leftarrow $ from user \label{alg_1_get_omR}
\Repeat  \label{alg_1_start_loop}
\State  $\theta,\pi\leftarrow $Solve \eqref{orig_RMP} over $\Omega_R$ \label{alg_1_RMP}
\State $l_* \leftarrow \mbox{arg}\min_{l \in \Omega}\bar{c}_l$ \label{alg_1_gen_col}
\State $\Omega_R \leftarrow \Omega_R \cup l_*$ \label{alg_1_add_col}
 \Until{$\bar{c}_{l_*} \geq 0$} \label{alg_1_end_loop}
 \State Return last $\theta$  generated \label{alg_1_ret}
\end{algorithmic}
\label{basicCG}
\end{algorithm} 

\begin{itemize}
    \item Line \ref{alg_1_get_omR}:  We initialize CG with columns from the user, which may consist of artificial variables with prohibitively high cost that ensure a feasible solution exists.  
    \item Line \ref{alg_1_start_loop}-\ref{alg_1_end_loop}:  We generate a sufficient set of columns to solve the MP exactly. 
    \begin{itemize}
        \item Line \ref{alg_1_RMP}:  Produce the solution to the RMP in \eqref{orig_RMP}.
        \item Line \ref{alg_1_gen_col}:  Produce the lowest reduced cost column.  CG does not require that \eqref{pricing} is solved exactly.  As long as a negative reduced cost column is added during any round of pricing (if such a column exists) then CG is guaranteed to solve \eqref{orig_LP} exactly \citep{costa2019,wang2017tracking}.
        \item Line \ref{alg_1_add_col}:  We add the generated column to $\Omega_R$. 
        Often more than one column is added during pricing.  This is facilitated in cases where pricing is solved using a dynamic program; as dynamic programs generate many solutions over the course of optimization.
        \item Line \ref{alg_1_end_loop}:  If no negative reduced cost column exists we terminate optimization.
    \end{itemize}
    \item Line \ref{alg_1_ret}:  We return the optimal solution.  This can be provided as input to a branch-price solver \citep{barnprice}, that calls Alg \ref{basicCG} in the inner loop.  
\end{itemize}

\section{Formal Description of Graph Generation}
\label{Sec_GG}
In this section we introduce our Graph Generation  (GG) algorithm. GG is an enhanced Column Generation (CG) algorithm that differs from CG by solving a more computationally intensive restricted master problem (RMP) at each iteration.  GG provides dual solutions that enforce many constraints that must be satisfied at optimality but have not yet been explicitly generated during pricing. 
This accelerates the convergence of CG as fewer iterations of CG are required. GG does not alter the structure of the pricing problem; nor does it loosen the master problem.
\subsection{Families of Columns}
\label{Sec_GG_fam}
In this subsection we describe the concept of families of columns, which we use in the remainder of this section.  Let $F$ be a set of subsets of $\Omega$, each member of which is called a family of columns.  We index $F$ using $f$ where $\Omega_f \subseteq \Omega$ is the family associated with $f$.  For any $l \in \Omega$ let $f_l$ be the family corresponding to $l$.  For any $l \in \Omega$ the associated family includes $l$ (meaning $l \in \Omega_{f_l} \quad \forall l \in \Omega$).

Given any $\Omega_R \subseteq \Omega$, let $\Omega_{R2} $ be the union of columns in the families of columns in $\Omega_R$; meaning $\Omega_{R2} = \cup_{l \in \Omega_R}\Omega_{f_l}$. With the aim of accelerating the convergence of CG we seek to solve $\Psi(\Omega_{R2})$ efficiently at each iteration of CG.  If $|\Omega_{R2}|$ is much larger than $|\Omega_R|$ then trivial methods would not be easily able to enumerate $\Omega_{R2}$ much less solve $\Psi(\Omega_{R2})$.  If solving for $\Psi(\Omega_{R2})$ is feasible it may or may not be the case that $\Psi(\Omega_R)-\Psi(\Omega_{R2})$ is large. However if $\Psi(\Omega_R)-\Psi(\Omega_{R2})$ is large over the course of CG and $\Psi(\Omega_{R2})$ can be solved easily at each iteration, then we hypothesize that CG converges faster than normal.  We validate this experimentally in this paper.  
%
\subsection{Path Cone}
\label{Sec_GG_pathCOne}
In order for solving $\Psi(\Omega_{R2})$ over \eqref{orig_RMP} to lead to faster convergence of CG in terms of time (not merely iterations) it must be the case that it is not an undue burden to solve $\Psi(\Omega_{R2})$. We now provide additional notation permitting us to express an LP formulation for fast solution of $\Psi(\Omega_{R2})$.  Consider that for any $l \in \Omega$ we can map it to an $f\in F$ s.t. $f=f_l$ for which the following terms are easily provided.
\begin{itemize}
    \item Let graph $G^f$ be associated with edge set $E^f$ and vertex set $V^f$. We index $E^f$ with $i,j$ corresponding to vertices in $V^f$.  There are special vertices  $v^+$ and $v^-$ in $G^f$ called the source and sink respectively.  The graph $G^f$ is directed and acyclic.  
    \item Let $P^f$ be the set of paths in $G^f$ starting at $v^+$ and ending at $v^-$, which we index by $p$.  There is a surjection from $P^f$ to $\Omega_f$, where $p$ maps to $l_p$.  We set $a_{ijp}=1$ if path $p$ uses edge $ij$ and otherwise set $a_{ijp}=0$.  A path $p$ satisfies the following flow constraint written using $[]$ to denote the binary indicator function. 
\begin{align}
    \sum_{ij\in E^f}a_{ijp}-\sum_{ji \in E^f}a_{jip}=[i=v_+]-[i=v_-] \quad \forall i \in V^f, p \in P^f, f \in F
\end{align}
We now express the convex cone of paths, which we refer to as path cone using non-negative values $\{ \psi^f_{ij} \quad \forall ij \in E^f\}$.  The path cone is the set of possible settings of vector $\psi^f$ (over all $ij \in E^f$) corresponding to a non-negative combination of paths.  Any non-negative vector $\psi^f$  lies in the path cone if there exists a non-negative vector of $\{\alpha_p \quad  \forall p \in P^f\}$ satisfying the following.     
\begin{align}
    \psi^f_{ij}=\sum_{p \in P^f}a_{ijp}\alpha_p \quad \forall ij \in E^f,f \in F
\end{align}
The path cone can be alternatively expressed as follows by relying on path, flow equivalence. 
\begin{align}
    \sum_{ij\in E^f}\psi^f_{ij}-\sum_{ji \in E^f}\psi^f_{ij}=0 \quad \forall i \in V^f-(v_+ \cup v_-), f \in F
\end{align}
    \item Each edge in $ij \in E^f$; is equipped with a cost $c^f_{ij}$ s.t. for any $p \in P^f$ the total cost of edges on the path $p$ is identical to that of $c_{l_p}$.  We write this formally below.
    \begin{align}
        c_{l_p}=\sum_{ij \in E^f}c^f_{ij}a_{ijp} \quad \forall p \in P^f, f \in F
    \end{align}
    \item Each edge $ij \in E^f$ is equipped with a vector $\vec{h}^f_{ij}$ with number of rows equal to the number of rows of $A$ s.t. the following holds.
    \begin{align}
        A_{:l_p}=\sum_{ij \in E^f}\vec{h}^f_{ij}a_{ijp} \quad \forall p \in P_f
    \end{align}
    Thus we can map a non-negative vector $\psi^f$ that lies in the path cone to the corresponding contribution to the constraints of the RMP.  
    \end{itemize}
%
\subsection{Graph Generation Algorithm}
\label{sec_GG_CGSOlve}
In this section we describe our Graph Generation (GG) algorithm for solving  expanded LP relaxations efficiently. Using our definitions in Section \ref{Sec_GG_pathCOne} we write $\Psi(\Omega_{R2})$ below over the path cone using $F_R=\cup_{l \in \Omega_R}f_l$ to be the union of families of columns in $\Omega_R$.
\begin{subequations}
\label{aug_RMP3}
\begin{align}
\Psi(\Omega_{R2})=    \min_{\substack{\theta \geq 0\\ \psi \geq 0 }}\sum_{l \in \Omega_R} c_l \theta_l+\sum_{\substack{f \in F_R\\ij \in E^f}}c^f_{ij} \psi^f_{ij}\\
    \sum_{l \in \Omega_{R}}A_{:l}\theta_l+\sum_{\substack{f \in F_R\\ij \in E^f}}\vec{h}^f_{ij}\psi^f_{ij} \geq \vec{b} \quad [\pi]\\
    \sum_{\substack{j\in V^f\\  ij \in E^f}} \psi^f_{ij}=\sum_{\substack{j\in V^f\\  ji \in E^f}} \psi^f_{ij} \quad \forall i \in V^f-(v_+,v_-), f \in F_R
\end{align}
\end{subequations}
Solving \eqref{orig_LP} is accomplished via CG using \eqref{aug_RMP3} to provide a primal/dual solution at each iteration.  Here pricing is never done on $\psi$ since $F_R$ is fixed as $F_R=\cup_{l \in \Omega_R}f_l$, which grows with $\Omega_R$.  Pricing over $\Omega_R$ is done using \eqref{pricing}.  
In Alg \ref{advCG} we write the CG optimization, which we refer to as the Graph Generation algorithm (GG), and provide annotation below.
\begin{algorithm}[!b]
 \caption{Graph Generation Algorithm (GG)}
\begin{algorithmic}[1] 
\State $\Omega_R\leftarrow $ from user \label{alg_2_init_1}
\State $F_R \leftarrow \Omega_R$ \label{alg_2_init_2}
\Repeat \label{alg_2_loop_start}
\State   $\theta,\psi,\pi \leftarrow$ Solve $ \Psi(\Omega_{R2})$ via \eqref{aug_RMP3}  \label{alg_2_solve_RMP}
\State $l_* \leftarrow \mbox{arg}\min_{l \in \Omega}\bar{c}_l$ \label{alg_2_pricing}
\State $\Omega_R \leftarrow \Omega_R \cup l_*$ \label{alg_2_add_omega}
\State $F_R \leftarrow F_R \cup f_{l_*}$ \label{alg_2_add_F}
 \Until{$\bar{c}_{l_*} \geq 0$} \label{alg_2_loop_end}
 \State Return last $\theta,\psi$  generated.  \label{alg_2_return}
\end{algorithmic}
\label{advCG}
\end{algorithm}

\begin{itemize}
    \item Line \ref{alg_2_init_1}-\ref{alg_2_init_2}: We initialize CG with columns from the user, which may consist of artificial variables with prohibitively high cost that ensure a feasible solution exists.  
    We initialize $F_R$ to $\Omega_R$.  If artificial primal variables (columns) are provided in $\Omega_R$  then these are not added to $F_R$.
    \item Line \ref{alg_2_loop_start}-\ref{alg_2_loop_end}:  We generate a sufficient set of columns/graphs to solve the MP exactly.  
    \begin{itemize}
        \item Line \ref{alg_2_solve_RMP}:  Solve the RMP providing a primal/dual solution.
        \item Line \ref{alg_2_pricing}:  Call pricing to generate the lowest reduced cost column $l_*$.  As in basic column generation we can generate more than one column.
    \item Line \ref{alg_2_add_omega}-\ref{alg_2_add_F}: Add the new column $l_*$ and the associated family $f_{l_*}$ to the RMP.
    \item Line \ref{alg_2_loop_end}:  If no column has negative reduced cost then we terminate optimization.      \end{itemize}

    \item Line \ref{alg_2_return}: Return last primal solution generated which is optimal.  As in basic CG this can be returned to a branch-price algorithm, that uses Alg \ref{advCG} as an inner loop operation.   
\end{itemize}
%

\section{Graph Generation for Capacitated Vehicle Routing}
\label{Sec_CVRP_apply}

In this section we apply Graph Generation (GG) to the Capacitated Vehicle Routing Problem (CVRP). We organize this section as follows.  In Section \ref{CVRP_Sec_Inf_desc} we provide a description of the CVRP problem. In Section \ref{CVRP_Sec_formal_desc} we provide a formal description of CVRP along with its master problem. In Section \ref{pricingILP} we describe pricing as an integer linear program, which may or may not be solved as such, and is typically solved as a resource constrained shortest path problem using a dynamic programming based labeling algorithm \citep{costa2019}. 
In Section \ref{CVRP_Sec_top_fam} we describe expanded topological families, which are the families that we use for GG in our application.  In Section \ref{CVRP_Sec_dt} we describe the application specific mechanism to generate the family of a column generated during pricing.  


\subsection{CVRP Problem Description}
\label{CVRP_Sec_Inf_desc}
 We are given  set of customers that must be serviced each with a location and demand; the location of a depot where vehicles start and end; and a set of homogeneous vehicles. We seek to assign vehicles to routes so as to minimize the total distance traveled; such that all customers are serviced and the capacity of the vehicles is respected.  


\subsection{Formal Description}
\label{CVRP_Sec_formal_desc}

%
We define the set of customers as $N$, which we index by $u$.  We use $N^+$ to denote $N$ augmented with the depot.  We have access to $K$ homogeneous vehicles each with capacity $d_0 \in \mathbb{Z}_+$, which start (and end) at the depot.  The demand of customer $u$ is denoted $d_u \in \mathbb{Z}_{+}$. We use $c_{uv}$ to denote the distance between any pair $u,v$ each of which lie in $N^+$.   
We use $\Omega$ to denote the set of feasible routes, which we index by $l$. We describe $l$ using the following notation.  
\begin{itemize}
    \item We set $a_{ul}=1$ if route $l$ services customer $u$, and otherwise set $a_{ul}=0$ for any $u \in N$.  
    \item We set $a_{uvl}=1$ if $u$ is immediately proceeded by $u$ in the route $l$, and otherwise set $a_{uvl}=0$ for any $u\in N^+,v \in N^+$. 
    \item 
    Below we define $c_l$ to be the cost of route $l$ where $c_l$ is defined as the total travel distance on route $l$.
    \begin{align}
        c_l=\sum_{\substack{u \in N^+\\ v \in N^+}}a_{uvl}c_{uv} \quad \forall l \in \Omega
    \end{align}
\end{itemize}
We write the master problem (MP) for CVRP formally below (with exposition provided below the equations).
\begin{subequations}
\label{RMP_CVRPTW}
\begin{align}
\min_{\theta \geq 0}\sum_{l \in \Omega}c_l\theta_l \label{eq_rmp_cvrp_obj}\\
\sum_{l \in \Omega}a_{ul}\theta_l \geq 1 \quad \forall u \in N  \label{eq_rmp_cvrp_cover}\\
\sum_{l \in \Omega}-\theta_l\geq -K  \label{eq_rmp_cvrp_pack}
\end{align}
\end{subequations}

In \eqref{eq_rmp_cvrp_obj} we seek to minimize the total distance traveled.  In \eqref{eq_rmp_cvrp_cover} we enforce that each customer is serviced at least once.  In \eqref{eq_rmp_cvrp_pack} we enforce that no more than $K$ vehicles are used where $K$ is the user defined number of vehicles.  Since $|\Omega|$ grows exponentially in the number of customers, column generation (CG) based methods are often applied to solve \eqref{RMP_CVRPTW}.

\subsection{Pricing as an Integer Linear Program}
\label{pricingILP}
We now consider the solution to pricing for CVRP as an integer linear program (ILP).  This section is adapted from \citep{yarkony_Detour_DOI} (Appendix B) for the convenience for the reader.  We use decision variable $x_{uvd}=1$ if the generated route services $u$ then goes to $v$ and contains (has remaining) exactly $d$ units of capacity after leaving $u$ and otherwise set $x_{uvd}=0$. The following combinations of $u,v,d$ exist.
\begin{itemize}
    \item $x_{uvd}$ exists if $d_v\leq d\leq d_0-d_u$ for any $u\neq v$ each of which lie in $N$.   
    \item $x_{-1,u,d_0}$ exists for each $u \in N$. This connects the depot to $u$ at the start of the route.  We use $-1$ to refer to the depot at the start of the route. Commas are added for convenience of the reader here.   Since we have only one depot, the start/end depot  are mathematical abstractions that refer to the  same place physically.  Also note that clearly $d_0\geq d_u$ as otherwise the problem would be infeasible.  
    \item $x_{u,-2,d}$ exists for each $u \in N,0\leq d \leq d_0-d_u$. This connects the $u$ to the depot at the end of the route.   We use $-2$ to refer to the depot at the end of the route. 
\end{itemize}
The set of valid combinations of $u,v,d$ is denoted $Q$.  
Below we define $\bar{c}_{uv}$ to be the cost of traveling from $u$ to $v$ minus the additional cost corresponding to dual variables so that  $\bar{c}_l=\sum_{(uvd)\in Q}\bar{c}_{uv}a_{uvl} ,\quad \forall l \in \Omega $ which facilitates the efficient writing of pricing.  
\begin{subequations}
\begin{align}
    \bar{c}_{uv}=c_{uv}-\pi_v \quad  \forall (uvd) \in Q, v \neq -2\\
    \bar{c}_{uv}=c_{uv}+\pi_0 \quad \forall (uvd) \in Q, v = -2 
\end{align}
\end{subequations}
We write the computation of the lowest reduced cost route below in the form of an ILP, which we annotate after the ILP.  
\begin{subequations}
\label{pricingEqFull}
\begin{align}
    \min_{x \in \{0,1\} }\sum_{\substack{(uvd)\in Q}}\bar{c}_{uv}x_{uvd} \label{ObjPricer}\\
    \sum_{\substack{(uvd) \in Q}}x_{uvd}= 1 \quad \forall u \in N \label{elemnt}\\
    \sum_{\substack{(-1,v,d) \in Q}}x_{-1vd}\leq 1 \label{pushCon}\\ 
    \sum_{\substack{(uvd) \in Q}}x_{uvd}=\sum_{\substack{(v,u,d-d_v) \in Q}}x_{v,u,d-d_v} \quad \forall  v\in N,d\geq d_v \label{flowcon}
\end{align}
\end{subequations}

In \eqref{ObjPricer} we minimize the reduced cost of the generated route.  In \eqref{elemnt} we ensure that a customer is visited no more than once. In \eqref{pushCon} we enforce that the vehicle leaves the start depot exactly once.  In \eqref{flowcon} we enforce that the vehicle leaves each customer it services with the appropriate amount of demand.    

We should note that the solution of \eqref{pricingEqFull} is not typically solved as an ILP but instead tackled as with labeling algorithm \citep{costa2019} for sake of efficiency.  The resources that must be kept track of by the labeling algorithm are \textbf{(a)} the set of customers visited thus far \textbf{(b)} the total amount of capacity used.  

\subsection{Expanded Topological Families}
\label{CVRP_Sec_top_fam}

We use an expanded variant of topological families inspired by \citep{haghani2021family}. Each family $f$ is associated with ordered list containing all $u \in N$. We describe the list associated with family $f$ using $\beta^f_u \in \mathbb{Z}_{+}$ where $\beta^f_u$ is the position in the ordered list that $u$ occupies.  A column $l$ lies in $\Omega_f$ if for any $u\in N,v\in N$ s.t.  $a_{uvl}=1$ then $u$ must come before $v$ in the ordering associated with family $f$.  Formally we write this as follows. 
\begin{align}
    \quad (l \in \Omega_f)\leftrightarrow \{ (a_{uvl}=1) \rightarrow (\beta^f_u <\beta^f_v) \quad \forall u \in N,v \in N\}\quad \forall f \in F, l \in \Omega \}
\end{align}
%
%
%
%
We describe the graph $G^f$ with $V^f,E^f$ used in GG as follows.  For every $u \in N$, $d \in [0,1,2...,d_0-d_u]$ create one vertex in $V^f$ denoted $(u,d)$.  We construct $E^f$ and the corresponding cost terms $c_{ij}^f$ as follows.  Connect vertex $(u,d)$ to vertex $(v,d-d_v)$ with an edge of cost $c_{uv}$ IFF $\beta^l_u<\beta^l_v$ (and $d-d_v\geq 0$). Connect every vertex $(u,d)$ to the sink $v^-$ with an edge of cost $c_{u,-2}$; where $c_{u,-2}$ is the distance from $u$ to the depot.  For every $u \in N$ connect the source vertex $v^+$ to vertex $(u,d_0-d_u)$ with an edge of cost $c_{-1,u}$;  where $c_{-1,u}$ is the distance from the depot to $u$. 

 We now consider the construction of $\{\vec{h}^f_{ij} \quad \forall ij \in E^f\}$.  We index the constraints in the MP for CVRP in \eqref{RMP_CVRPTW} with $u\in N$ and $0$ for the the constraint enforcing no more than $K$ routes used.  We denote the corresponding entries in $\vec{h}^f_{ij}$ as $\vec{h}^f_{ij;u} $ and $\vec{h}^f_{ij;0}$ respectively.  For any $\vec{h}^f_{ij}$ where $j$ corresponds to $(u,d)$ for some $(u,d)$ we define $\vec{h}^f_{ij;u}=1$.  All other entries of $\vec{h}^f_{ij}$ are zero for that $ij \in E^f$.  For $j=v^-$ we define $\vec{h}^f_{ij;0}=-1$. All other entries of $\vec{h}^f_{ij}$ are zero for that $ij \in E^f$.

\subsection{Heuristic Construction of the Expanded Topological Family}
\label{CVRP_Sec_dt}
 Given any $l \in \Omega$; in this section we consider the heuristic construction of the ordering described using $\beta^{f_l}$ so that $\Omega_{f_l}$ contains useful routes.  Observe that items that are in similar physical locations should be in similar positions on the ordered list so that a route in $\Omega_{f_l}$ can visit all of them without leaving the area and then coming back. We initialize the ordering with the items in $N_l$ sorted in order from first visited to last visited.  This way $l$ lies in $\Omega_{f_l}$.  
 Now iterate over the remaining customers ($N-N_l$) in a random order.  For each customer $u$ insert it immediately behind the customer in $N_l$ that is closest to $u$. For items closer to the depot than any customer in $N_l$ insert it in the beginning of the list.  Observe that after each insertion $l$ lies in $\Omega_{f_l}$.  

\section{Experimental Validation of Graph Generation}
\label{Sec_exper}
In this section we demonstrate the effectiveness of our Graph Generation (GG) for the Capacitated Vehicle Routing Problem (CVRP). We compared the performance of GG against standard  (meaning unstabilized) Column Generation (CG).  


Following the example of \citep{yarkony_Detour_DOI}, to provide fair comparisons, a ``vanilla" implementation of pricing is implemented.  We generate one column at each iteration of CG (or GG); where this column is the lowest reduced cost column.  Hence no heuristics are used during pricing.  To do this we solve the integer linear programming (ILP) formulation of pricing described in \citep{yarkony_Detour_DOI},  as described in Section \ref{pricingILP}.

For our experiments we considered 25 randomly generated problem instances of the following form.  Each problem instance is associated with 30 customers of demand one and 5 vehicles of capacity 7 each. Each customer and the depot is assigned a random integer position on grid of size 100 by 100. Distances between customers (and  the depot) are computed based on the L2 (Euclidian) distance rounding up to the nearest integer. 

To solve the restricted master problem (RMP) during the course of CG optimization we used the basic MATLAB linear programming solver with default options. For pricing we used the intlinprog in MATLAB with default options.  We initialize CG with one artificial column for each item. This column uses no vehicle; has prohibitively high cost; and covers the corresponding item.  
%
%
%

\begin{figure}[!htpb]
\centering
	\includegraphics[width=0.49\linewidth]{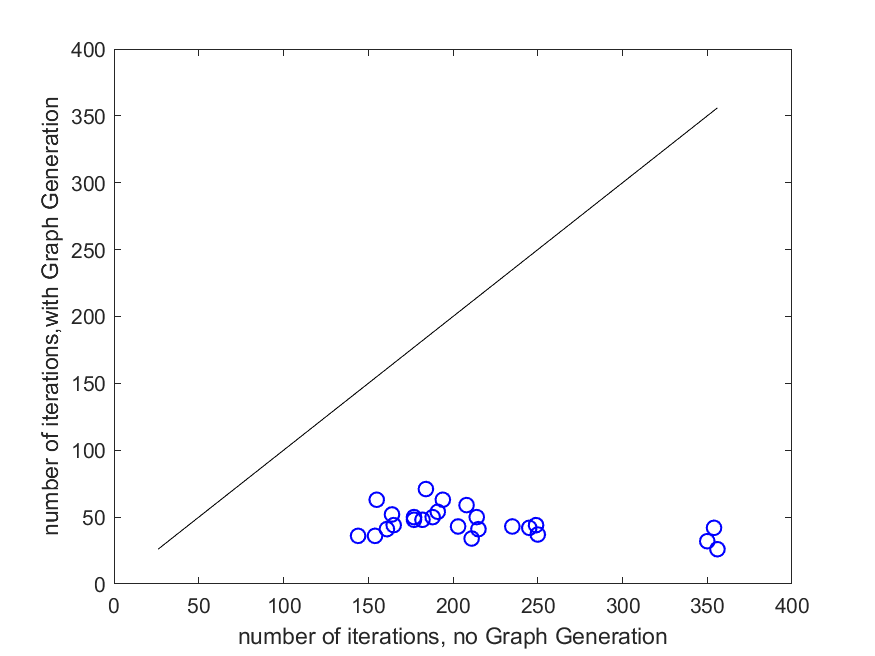}
	\includegraphics[width=0.49\linewidth]{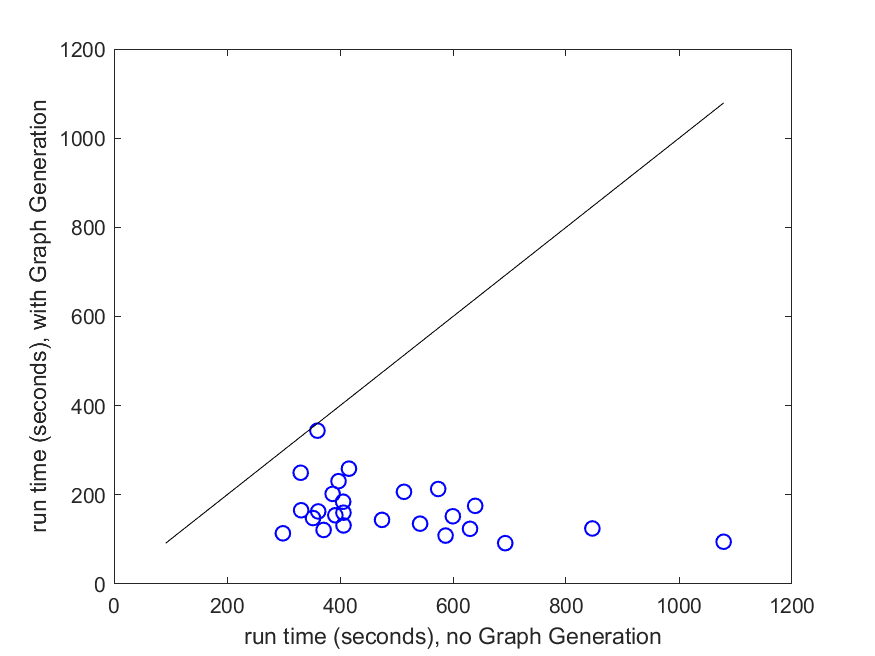}
	\caption[Detour-DOI results]{
		\textbf{(Left):  }Comparative iterations required between using GG vs no stabilization for all problem instances. Each blue dot describes the performance of CG and GG on a single problem instance.  Its x-coordinate is the number of iterations CG took to solve that problem instance; while the y-coordinate is the number of iterations required to solve that problem instance with GG.  The black line plots the line y=x so as to provide a baseline for improvement.  The further below the line a point is the greater the improvement achieved by GG over CG.  
		\textbf{(Right):} Comparative run times between using GG and CG.  Lines and dots describe the same terms as (Left) except considering run time.}
	\label{fig:CVRP_results}
\end{figure}

\begin{figure}[!hbtp]
	\includegraphics[width=0.49\linewidth]{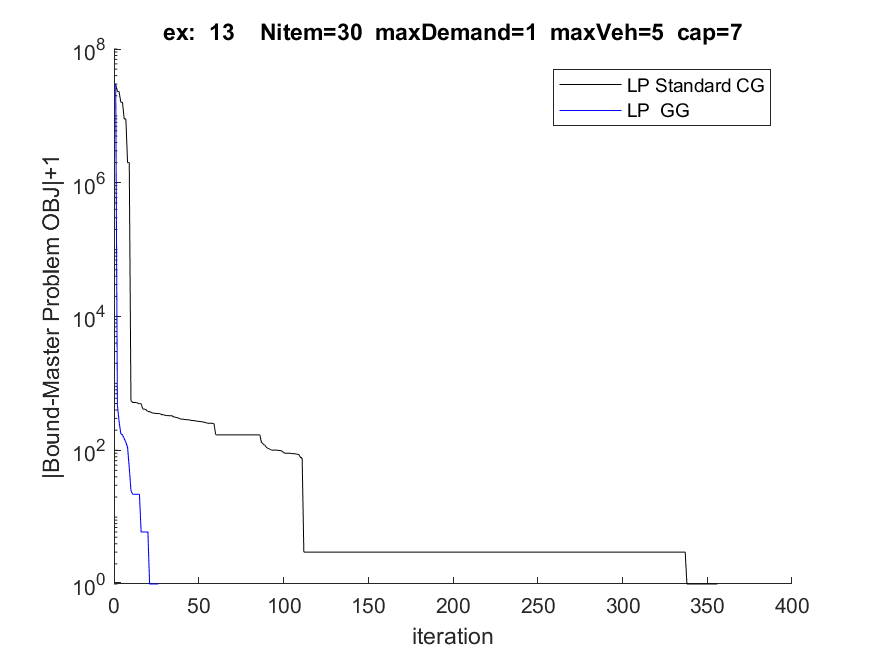}
	\includegraphics[width=0.49\linewidth]{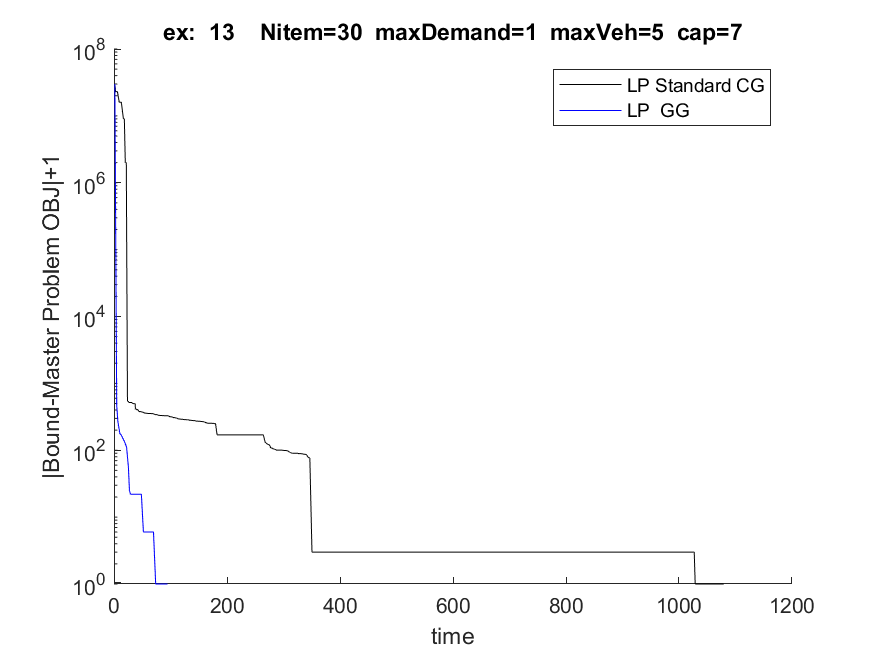}\\
		\includegraphics[width=0.49\linewidth]{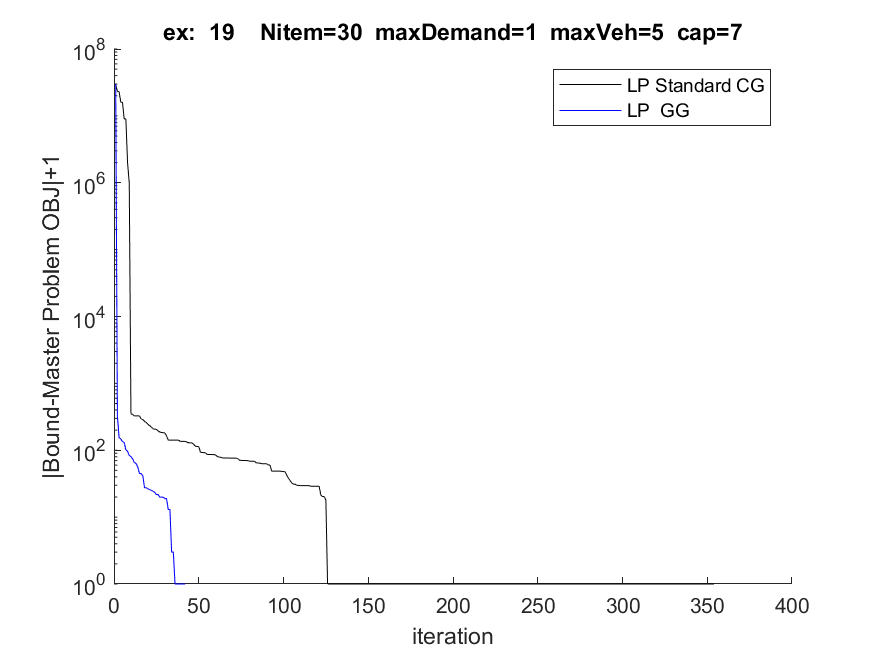}
	\includegraphics[width=0.49\linewidth]{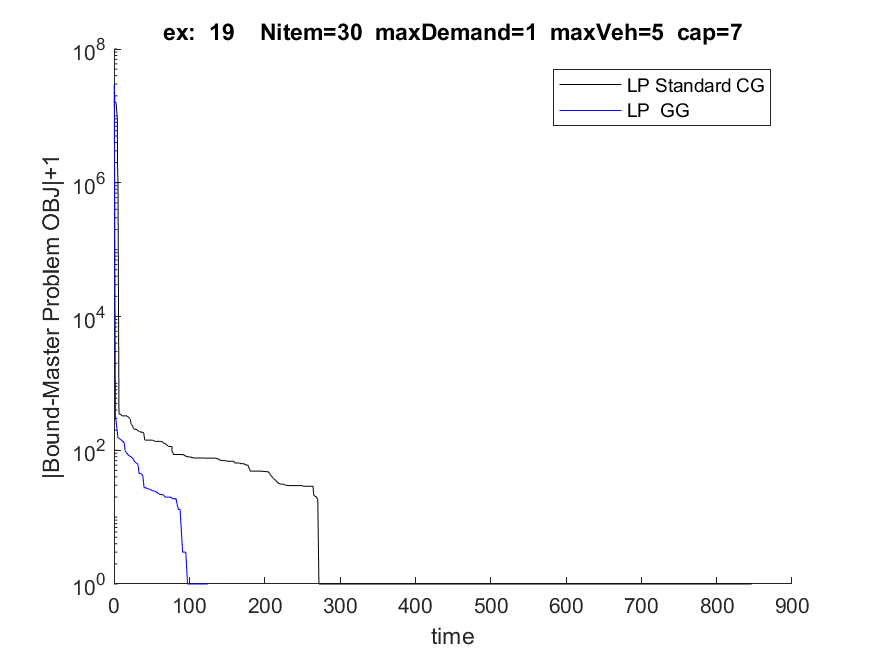}\\
		\includegraphics[width=0.49\linewidth]{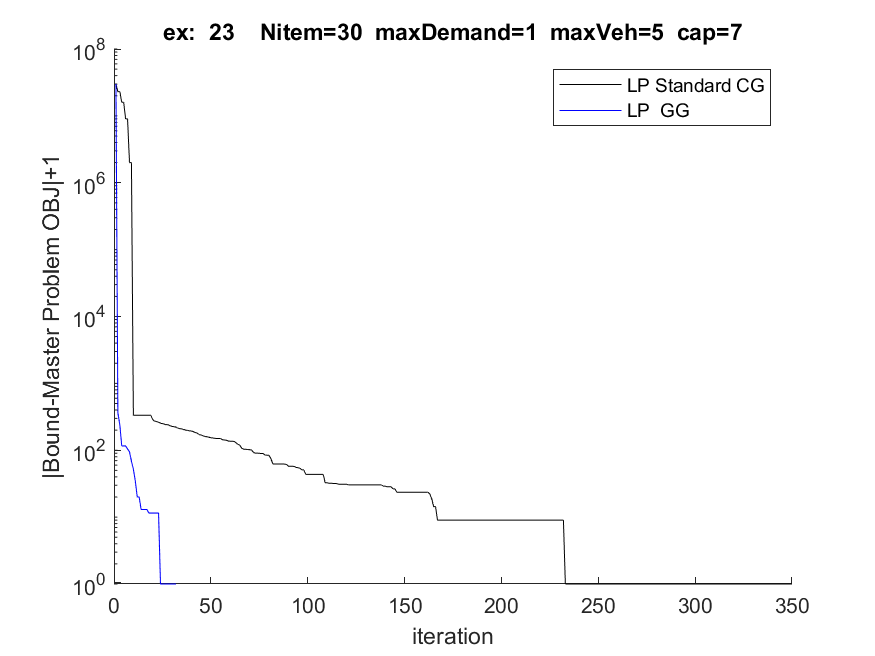}
	\includegraphics[width=0.49\linewidth]{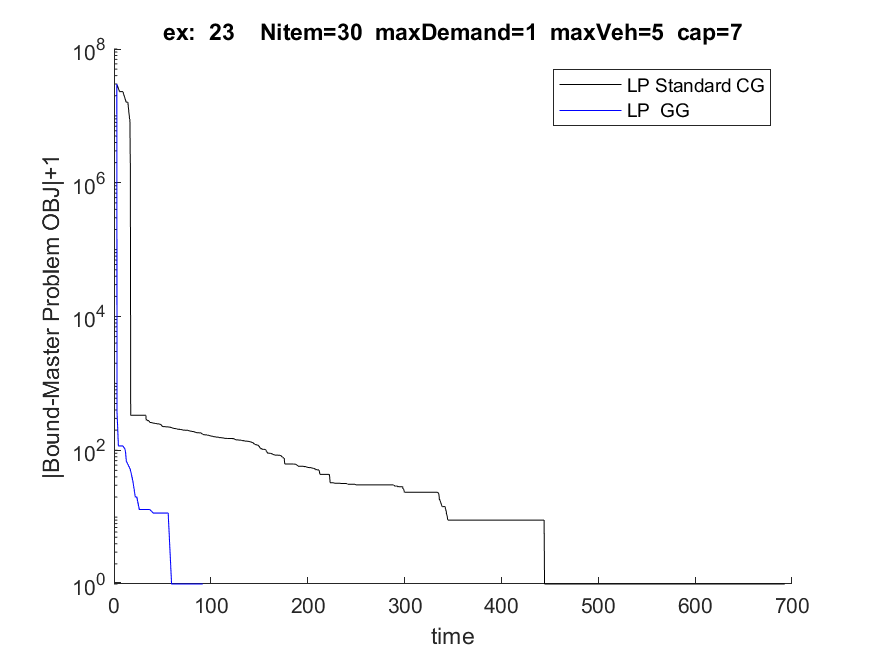}\\
	\caption{Results on individual problem instances for the value of the RMP LP as a function of iteration/time (sec) in semilog scale. The left side provides the results with respect to iteration and the right side for time on the same instance. We add one to the difference the LP RMP and the MP values, which allows us to use the semilog scale.
	}
	\label{fig:iter2}
\end{figure}

In Tables  \ref{table:CVRP_results_iters} and \ref{table:CVRP_results_time} we provide the number of iterations and time required to achieve the optimal master problem (MP) solution for CG and GG on our problem instances.  This data is aggregated in plots in Fig \ref{fig:CVRP_results}.  We observe that GG always outperforms CG in terms of both time and iterations and that  the amount of performance improvement by GG over CG grows as the number of iterations/time taken by CG increases. In Fig \ref{fig:iter2} visualize convergence on sample problem instances displaying the value of the RMP LP as a function iteration/time.

\begin{table}[tpb]
	\centering
	\scalebox{0.85}{
	\begin{tabular}{|c|c c|} 
		\hline
		& \multicolumn{2}{c|}{\bf Iterations}\\
		\bf Instance & \bf un-stabilized CG & \bf Graph Generation \\
		\hline
        1 & 165 & 44 \\
        2 & 191 & 54 \\
        3 & 154 & 36 \\
        4 & 211 & 34 \\
        5 & 250 & 37 \\
        6 & 249 & 44 \\
        7 & 161 & 41 \\
        8 & 184 & 71 \\
        9 & 203 & 43 \\
        10 & 235 & 43 \\
        11 & 188 & 50 \\
        12 & 245 & 42 \\
        13 & 356 & 26 \\
        14 & 164 & 52 \\
        15 & 214 & 50 \\
        16 & 155 & 63 \\
        17 & 177 & 50 \\
        18 & 208 & 59 \\
        19 & 354 & 42 \\
        20 & 177 & 48 \\
        21 & 194 & 63 \\
        22 & 144 & 36 \\
        23 & 350 & 32 \\
        24 & 182 & 48 \\
        25 & 215 & 41 \\
        \hline
        mean & 213.0 & 46.0 \\
        median & 194.0 & 44.0 \\
		\hline
	\end{tabular}}
	\caption{CVRP Iteration Results}
	\label{table:CVRP_results_iters}
\end{table}

\begin{table}[tpb]
	\centering
	\scalebox{0.85}{
	\begin{tabular}{|c|c c|} 
		\hline
		& \multicolumn{2}{c|}{\bf Time (seconds)} \\
		\bf Instance & \bf Un-Stabilized CG & \bf Graph Generation \\
		\hline
            1 & 391.7 & 154.1 \\
            2 & 405.5 & 184.6 \\
            3 & 298.9 & 113.7 \\
            4 & 630.1 & 123.8 \\
            5 & 586.9 & 108.6 \\
            6 & 541.8 & 135.3 \\
            7 & 406.1 & 159.9 \\
            8 & 360.1 & 344.0 \\
            9 & 474.4 & 143.8 \\
            10 & 406.2 & 131.3 \\
            11 & 513.3 & 206.8 \\
            12 & 599.7 & 152.0 \\
            13 & 1079.0 & 94.7 \\
            14 & 386.9 & 202.1 \\
            15 & 573.8 & 213.1 \\
            16 & 330.5 & 249.5 \\
            17 & 331.2 & 165.4 \\
            18 & 397.3 & 230.6 \\
            19 & 846.7 & 124.5 \\
            20 & 352.2 & 147.9 \\
            21 & 415.8 & 258.7 \\
            22 & 371.0 & 121.2 \\
            23 & 692.5 & 91.7 \\
            24 & 361.4 & 162.8 \\
            25 & 639.2 & 175.3 \\
        \hline
            mean & 495.7 & 167.8 \\
            median & 406.2 & 154.1 \\

		\hline
	\end{tabular} }
	\caption{CVRP Runtime Results}
	\label{table:CVRP_results_time}
\end{table}

\section{Conclusion}
\label{Sec_conc}

In this document we expand on the Column Generation algorithm (CG) to solve expanded linear programming relaxations, producing the Graph Generation algorithm (GG).  GG is an accelerated form of CG that operates by providing improved solutions to the restricted master problem (RMP) at each iteration of CG.  GG solves a more computationally intensive master problem but is able to decrease the number of iterations of pricing needed.  Hence GG works well in problem domains where pricing ,not solving the RMP, is the computational bottleneck.  GG is distinct from CG in the following manner.  At each step of pricing, GG adds additional variables/constraints to the RMP permitting the description of large numbers of columns including the column generated during pricing.  This is done in a manner that does not explode the size of the RMP.  This is in contrast to CG which simply adds one or more columns with negative reduced cost to the RMP.  

We apply GG to the Capacitated Vehicle Routing Problem and demonstrate that GG achieves large speedups over CG. In future work we seek to apply GG to problems with large numbers of resources for which restocking only happens on certain occasions and resource consumption happens uniformly across resources. Such problems specifically where tasks are ordered (in time), naturally produce compact graphs, expressing large numbers of columns for GG to use. Crew scheduling is an example of such a problem where days off/on follow specific structural requirements, but the schedule of any given day is independent given that structure. 

In future work we seek to create principled mechanisms for the generation of orderings of items for various applications where such families are used in GG. In addition we seek to apply and adapt column management techniques \citep{lubbecke2005selected} to decrease the size of the RMP when solving the GG RMP becomes a computational bottleneck.  
We are experimenting with a number of related techniques on larger problem instances. 
Specifically we are adapting ideas in ``column management" \citep{Desaulniers2005}, which reduce the number of columns in the CG RMP, to GG.  These techniques which we refer to as ``graph management" methods, remove subsets of the graphs (or entire graphs) from the RMP that have not been active in the primal RMP in recent iterations.  Similarly graph management approaches may add in graphs or subsets of graphs that have been removed previously but may be helpful again.   
\par To permit us to clearly explain the overall method of GG we have excluded these improvement methods from discussion in this paper. Finally, we are also working on heuristic methods to generate good integer solutions if upon termination of GG we have not already found one. We do this by solving the GG RMP as an integer linear program. We are also studying the integration of GG into branch-price algorithms \citep{barnprice}.  

\bibliographystyle{abbrvnat} 
\singlespacing
\bibliography{col_gen_bib}

\end{document}